\documentclass[12pt]{elsarticle}

\usepackage{amssymb}
\usepackage{amssymb}
\usepackage{amssymb}
\usepackage{amssymb,graphics,float}
\usepackage[centertags]{amsmath}
\usepackage{pstricks,pst-node,pst-tree,pst-plot}
\usepackage{graphicx,epsfig,subfigure}
\usepackage{pstricks,pst-node,pst-tree,pst-plot}

\usepackage{amsthm}

\usepackage[labelsep=endash]{caption}
\usepackage{graphicx}
\usepackage{placeins}
\usepackage{tikz}

\makeatletter
\def\ps@pprintTitle{%
 \let\@oddhead\@empty
 \let\@evenhead\@empty
 \def\@oddfoot{\centerline{\thepage}}%
 \let\@evenfoot\@oddfoot}
\makeatother

\begin{document}

\begin{frontmatter}

\title{Control of Parasitism in Variational Integrators for Degenerate Lagrangian Systems}

\author[label1]{ Farrukh Shehzad}
\ead{farrukhmath@gmail.com}

\author[label1]{ Yousaf Habib\corref{cor1}}
\cortext[cor1]{Corresponding author}
\ead{yhabib@cuilahore.edu.pk}

\author[label2,label3]{ Michael Kraus}
\ead{michael.kraus@ipp.mpg.de}

\author[label1]{ Zareen Akhtar}
\ead{zareenakhtar124@yahoo.com}

\address[label1]{Department of Mathematics, COMSATS University Islamabad, Lahore Campus, Pakistan.}

\address[label2]{Max-Planck-Institut für Plasmaphysik,
Boltzmannstraße 2, 85748 Garching, Deutschland.
}
\address[label3]{Technische Universität München, Zentrum für Mathematik,
Boltzmannstraße 3, 85748 Garching, Deutschland.
}
\date{}

\begin{abstract}
This paper deals with the control of parasitism in variational integrators for degenerate Lagrangian systems by writing them as general linear methods. This enables us to calculate their parasitic growth parameters which are responsible for the loss of long-time energy conservation properties of these algorithms. As a remedy and to offset the effects of parasitism, the standard projection technique is then applied to the general linear methods to numerically preserve the invariants of the degenerate Lagrangian systems by projecting the solution onto the desired manifold.

\end{abstract}
 

\end{frontmatter}


\section{Introduction}

\noindent Variational integrators discretise the action integral of the Lagrangian $L(q, \dot{q})$ of a dynamical system, where $q$ is position and $\dot{q}$ is velocity. A discrete analogue of Hamilton's principle of stationary action is then applied, which results in the discrete Euler Lagrange equations and the corresponding evolution map is termed as variational integrators \cite{NA16,NA17, NA18, NA19,11}.\\
In various problems of Physics, we deal with degenerate Lagrangian systems, whose Lagrangian is of the form 
\begin{equation}\label{dl}
     L= \left<\alpha(q),\dot{q}\right>-H(q) , 
\end{equation}
such that,
$$ \frac{\partial ^2 L}{\partial{\dot{q}^i}{\partial{\dot{q}}^j}} = 0,$$
where $\alpha(q)$ is possibly a non-linear function of $q$ and $H(q)$ is the Hamiltonian of the system. Examples of such systems include the non-linear pendulum, planar point vertices and guiding centre dynamics. Applying standard variational integrators to such a degenerate Lagrangian system leads to multistep numerical methods, which are prone to parasitic instabilities. The reason being that the integration process magnifies the perturbation in non-principle parasitic components of the numerical solution and leads to numerical corruption. Degenerate variational integrators \cite{2,1} allow to construct one-step methods for degenerate Lagrangians, but are limited to special forms of \eqref{dl} and not generally applicable.
In this paper, we write the resulting multistep methods as general linear methods and then apply the projection technique \cite{NA0,6,7} to control the effects of parasitism \cite{rfb18,rfh1}. This approach has the advantage that the resulting methods are applicable to all degenerate Lagrangians of the form \eqref{dl}.

\section{Variational integrators}
\begin{itemize}
\item \noindent The first step in the construction of variational integrators is to discretise the action integral given as,
\begin{equation}\label{act_int} A[q(t)]=\int{L(q(t), \dot{q}(t)}) dt.
\end{equation}
\noindent where,   $$ L: TQ \xrightarrow{}{}  R, $$ denotes the Lagrangian of a mechanical system with configuration manifold $Q$,
TQ is the tangent bundle and represents the velocity phase space, q(t) denotes the trajectory in Q and $\dot{q}(t)$ its time derivative. 

\item \noindent The second step is to find the discrete analogue to the action integral \eqref{act_int}. For this purpose, we divide the time interval into an equidistant monotonic sequence ${\{q_m\}}_{m=0}^{N}$ joined by 
 a discrete curve and then add the discrete Lagrangian $ L_d (q_m,q_{m+1})$ on each adjacent pair to obtain the discrete action given as, 
 $$ A_d[q_d]=\sum_{m=0}^{N-1} L_d (q_m,q_{m+1}) ,$$
 where the discrete Lagrangian is,
\begin{equation}\label{dis_int}
 L_d(q_m, q_{m+1} ) \approx  \int_{t_m}^{t_{m+1}} L(q(t),\dot{q}(t))dt.
\end{equation}

\noindent One way to obtain the discrete Lagrangian is to use finite differences to approximate the velocity as,
\begin{align*}
 \dot{q} \approx \frac{q_{m+1}-q_m}{h} .
\end{align*}

\noindent The integral in \eqref{dis_int} cannot be computed analytically, so we resort to numerical approximation using a quadrature rule such as the trapezoidal rule and obtain,
\begin{align} \label{NA2}
  L_d (q_{m},q_{m+1})& =\frac{h}{2}[ L(q_{m},\frac{q_{m+1}-q_{m}}{h}) + L(q_{m+1},\frac{q_{m+1}-q_{m}}{h})]
  \end{align}
where $ L_d: Q \times Q \xrightarrow{}{}  R$ .

\item \noindent The third step is to find the variation of the discrete action i.e.,
\begin{align*}
  \delta A_d[q_d] &=\delta\sum_{m=0}^{N-1} L_d (q_m,q_{m+1}), \\
   &=\sum_{m=0}^{N-1} [D_1 L_d (q_m,q_{m+1}) \cdot \delta{q_m}+D_2 L_d (q_m,q_{m+1}) \cdot \delta{q_{m+1}}].  
\end{align*}
\noindent Here we have used integration by parts, and $D_1$ and $D_2$ denote the derivative with respect to the first and the second arguments respectively.
\begin{align*}
  \delta A_d[q_d] &=D_1 L_d (q_0,q_1) \cdot \delta{q_0} + \sum_{m=1}^{N-1} D_1 L_d (q_m,q_{m+1})\cdot \delta{q_m}\\
  &+\sum_{m=0}^{N-2}D_2 L_d (q_m,q_{m+1}) \cdot \delta{q_{m+1}}+D_2 L_d (q_{N-1},q_{N}) \cdot \delta{q_{N-1}}.  
\end{align*}
\noindent Since $\delta{q_0}= \delta{q_N}=0$, because the variations at the endpoints are fixed, therefore,
\begin{align*}
  \delta A_d[q_d] &= \sum_{m=1}^{N-1} D_1 L_d (q_m,q_{m+1}).\delta{q_m}+\sum_{m=0}^{N-2}D_2 L_d (q_m,q_{m+1}).\delta{q_{m+1}}.  
  \\ &= \sum_{m=1}^{N-1}[ D_1 L_d (q_m,q_{m+1}).\delta{q_m}+D_2 L_d (q_{m-1},q_{m})].\delta{q_{m}}.  
\end{align*}
\item \noindent The fourth step is to apply Hamilton’s principle of stationary action which requires $\delta A_d = 0$ for all $\delta{q_{d}}, $ 
Thus we obtain discrete Euler-Lagrange equations,
\begin{equation} \label{N2}
 D_1 L_d (q_m,q_{m+1})+D_2 L_d (q_{m-1},q_{m})=0.
\end{equation}
\noindent which define an evolution map 
$$ \phi_h : Q \times Q \xrightarrow{}{}   Q \times Q : (q_{m-1},q_m ) \xrightarrow[]{} (q_{m},q_{m+1} ).$$ 
\noindent This map determines $q_{m+1}$ from the given values of $q_{0},q_{1},\hdots,q_{m}$.
\end{itemize}
\subsection{Variational integrators for degenerate Lagrangian systems}
\noindent For the special case of degenerate Lagrangian system \eqref{dl}, an application of variational integrator yields first order Euler Lagrange equations. Specifically, let us apply the trapezoidal rule discretisation \eqref{NA2} to the degenerate Lagrangian system \eqref{dl}, we have,
\begin{align*}
 L_d (q_{m},q_{m+1}) &=\frac{h}{2}[ \alpha(q_{m}) \cdot \frac{q_{m+1}-q_{m}}{h} -H(q_{m})+\alpha(q_{m+1})\cdot \frac{q_{m+1}-q_{m}}{h} -H(q_{m+1})],\\
  L_d (q_{m-1},q_{m}) &=\frac{h}{2}[ \alpha(q_{m-1})\cdot \frac{q_{m}-q_{m-1}}{h} -H(q_{m-1})+\alpha(q_{m})\cdot \frac{q_{m}-q_{m-1}}{h} -H(q_{m})].
\end{align*}

\noindent The differentiation yields, 
\begin{align*}
D_1 L_d (q_{m},q_{m+1}) &=\frac{h}{2}[ \nabla \alpha(q_{m})\cdot \frac{q_{m+1}-q_{m}}{h} -\frac{\alpha(q_{m})}{h}-\frac{\alpha(q_{m+1})}{h} -\nabla H(q_{m})],\\
D_2 L_d (q_{m-1},q_{m}) &=\frac{h}{2}[ \nabla \alpha(q_{m}).\frac{q_{m}-q_{m-1}}{h} + \frac{\alpha(q_{m})}{h}+\frac{\alpha (q_{m-1})}{h} -\nabla H(q_{m})].
\end{align*}

\noindent The discrete Euler-Lagrange equations \eqref{N2} thus become,
\begin{align}  \label{ABC2}
    \nabla \alpha (q_m).(q_{m+1}-q_{m-1})=\alpha(q_{m+1})-\alpha(q_{m-1})+2h\nabla H(q_m).
\end{align}
\noindent The variational integrator in equation \eqref{ABC2} represents a multistep method and hence suffers from parasitic instabilities. We aim to write it in the form of general linear methods and use projection techniques to counteract the effects of parasitism.
\section{General linear methods}
\noindent General linear methods are numerical methods to calculate approximate solutions of initial value problems \cite{rfb17,8},
\begin{equation}\label{basic_eq}
    y'=f(x,y), \hspace{0.5in} y(0)=y_0,
\end{equation}
where  $f:{\mathbb{R}}^{N}\rightarrow {\mathbb{R}}^N$ and $x\in \mathbb{R}$. The  general linear methods in its general form can be written as,
\begin{eqnarray}\label{GLM}
   K&=&h(A\otimes I)f(K)+(U\otimes I)y^{[m-1]},\\
    y^{[m]}&=&h(B\otimes I)f(K)+(V\otimes I)y^{[m-1]}\notag
\end{eqnarray}
with $s$-number of stages
 $K\in ({\mathbb{R}}^N)^s$ and $r$-component input vector $y^{[m-1]}\in ({\mathbb{R}}^N)^r$ and output vector $y^{[m]}\in ({\mathbb{R}}^N)^r$ given as,
\begin{align*}
    K &= \begin{bmatrix}
           K_{1} \\
           K_{2} \\
           \vdots \\
           K_{s}
         \end{bmatrix},
  \,\,\,\,\ 
  f(Y)=\begin{bmatrix}
           f(K_1) \\
           f(K_2) \\
           \vdots \\
           f(K_s)
         \end{bmatrix},
          \,\,\,\,\ 
  y^{[m-1]}= \begin{bmatrix}
           {y_1}^{[m-1]} \\
           {y_2}^{[m-1]} \\
           \vdots \\
           {y_r}^{[m-1]}
         \end{bmatrix},
          \,\,\,\,\ 
  y^{[m]}= \begin{bmatrix}
           {y_1}^{[m]} \\
           {y_2}^{[m]} \\
           \vdots \\
           {y_r}^{[m]}
         \end{bmatrix}
  \end{align*}
The characteristic matrices $(A,U,B,V)$ of a GLM are referred to as,
\begin{equation}\label{matrix_glm}
\left[ {\begin{array}{c|c}
 A & U \\ \hline
B & V 
 \end{array} } \right].
\end{equation}
General linear methods include Runge--Kutta methods and all multistep methods. The Runge--Kutta methods  have a single input with $r=1$ so that the matrices $U=\mathbf{1}$, $V=1$ and  $B$ has a single row. An example of a two stage Runge-Kutta method in general linear formulation is,
\[  \left[ {\begin{array}{c}
Y_{1} \\
Y_{2}\\ \hline
y^{[m]}
 \end{array} } \right]   = \left[ {\begin{array}{cc|c}
a_{11} & a_{12} & 1 \\
a_{21} & a_{22} & 1 \\ \hline
b_{1} & b_{2} & 1 \\  \end{array} } \right]\left[ {\begin{array}{c}
hf(Y_{1}) \\
hf(Y_{2}) \\\hline
y^{[m-1]}
 \end{array} } \right].
\]
The linear multistep methods such as Adams-Moulton method given as,
\begin{equation*}
y_{m} = y_{m-1} + h(\beta_{0}f(y_{m}) + \beta_{1}f(y_{m-1}) + \beta_{2}f(y_{m-2}) +\cdots+ \beta_{k}f(y_{m-k})),
\end{equation*} 
written in general linear method formulation has $s=1$ and is given as,
\begin{equation}\label{Adams_Moulton}
\left[ {\begin{array}{c}
Y_{1} \\ \hline
y_{m}\\
hf(Y_{1})\\
hf(y_{m-1})\\
hf(y_{m-2})\\
\vdots\\
hf(y_{m-k+1})\\
 \end{array} } \right]   = \left[ {\begin{array}{c|cccccc}
\beta_{0} & 1 &\beta_{1} & \beta_{2} & \cdots & \beta_{k-1} & \beta_{k}\\ \hline
\beta_{0} & 1 &\beta_{1} & \beta_{2} & \cdots & \beta_{k-1} & \beta_{k}\\
1 & 0 & 0& 0& \cdots& 0& 0\\ 
0 & 0 & 1& 0& \cdots& 0& 0\\ 
0 & 0 & 0& 1& \cdots& 0& 0\\ 
\vdots & \vdots & \vdots & \vdots & \ddots & \vdots & \vdots\\
0 & 0 & 0& 0& \cdots& 1& 0\\ 
 \end{array} }  \right]  \left[ {\begin{array}{c}
hf(Y_{1}) \\ \hline
y_{m-1}\\
hf(y_{m-1})\\
hf(y_{m-2})\\
hf(y_{m-3})\\
\vdots\\
hf(y_{m-k})\\
 \end{array} } \right] 
\end{equation} 

\section{Parasitism in general linear methods}
\noindent General linear methods suffer from parasitic solutions which are obtained in addition to the numerical approximation
of the exact solution. The main reason  is that the perturbation in parasitic  components of the numerical solution is amplified with the passage of time \cite{rfb18,rfv12,rfh7}. Let us consider a GLM,
\begin{align*}
    \begin{bmatrix}
     \begin{array}{c}
        Y_1  \\
         Y_2\\
         \hline
         y_1^{[m]}  \\
         y_2^{[m]}\\
     \end{array}
    \end{bmatrix}=\begin{bmatrix}
     \begin{array}{cc|cc}
        a_{11} & a_{12} & u_{11} &u_{12} \\
         a_{21} &a_{22} &u_{21} &u_{22} \\
         \hline
         b_{11} &b_{12} & v & 0   \\
         b_{21} &b_{22} & 0 & -v \\
     \end{array}
    \end{bmatrix}\begin{bmatrix}
     \begin{array}{c}
        hf(Y_1)  \\
        hf(Y_2)\\
         \hline
         y_1^{[m-1]}  \\
         y_2^{[m-2]}\\
     \end{array}
    \end{bmatrix}
\end{align*}
with $v=1$ are eigen values of $V$. Here $ y_1^{[m]}$ approximates actual solution and  $ y_2^{[m]}$ is the parasitic numerical solution.
\begin{align*}
    Y_1&=ha_{11}f(Y_1)+ha_{12}f(Y_2)+u_{11}y_1^{[m-1]}+u_{12}y_2^{[m-1]}\\
    Y_2&=ha_{21}f(Y_1)+ha_{22}f(Y_2)+u_{21}y_1^{[m-1]}+u_{22}y_2^{[m-1]}\\
    y_1^{[m]}&=hb_{11}f(Y_1)+hb_{12}f(Y_2)+y_1^{[m-1]}\\
    y_2^{[m]}&=hb_{21}f(Y_1)+hb_{22}f(Y_2)-y_2^{[m-1]}
\end{align*}
An induced perturbation in the parasitic component of the numerical solution gives,
\begin{align*}
    y_2^{[m-1]}\mapsto y_2^{[m-1]}+(-1)^{m-1}z_{m-1}
\end{align*}
This perturbation in the stages $Y_{i}$ and stage derivatives  $F_{i}$ are, 
\begin{align*}
\delta Y_{i} &= (-1)^{n-1}u_{i2}z_{m-1} ,\\
\delta F_{i} &= (-1)^{n-1}\frac{\partial f}{\partial y}u_{i2}z_{m-1} .
\end{align*}
And the perturbation in the second output component is,
\begin{align*}
y_{2}^{[m]} + (-1)^{n}z_{m} &=h\displaystyle\sum_{i=1}^{2}b_{2i}F_{i}- y_{2}^{[m-1]} + (-1)^{n-1}h\displaystyle\sum_{i=1}^{2}\frac{\partial f}{\partial y}b_{2i}u_{i2}z_{m-1} - (-1)^{n-1}z_{m-1},\\
\Rightarrow z_{m}&=\big(1-h\displaystyle\sum_{i=1}^{2}\frac{\partial f}{\partial y}b_{2i}u_{i2}\big)z_{m-1}.
\end{align*}
This is similar to Euler method for the solution of the differential equation,
\begin{align*}\label{para_euler}
z'=\mu \frac{\partial f}{\partial y} z,
\end{align*}
where $\mu = -\displaystyle\sum_{i=1}^{2}b_{2i}u_{i2} $ is the first order parasitic growth parameter of the general linear method and can be computed as,
\[
BU= \left[ \begin{array}{cc}
1 & 0  \\
0 & -\mu \end{array} \right].
\]
The second order parasitic growth parameter has been calculated in \cite{rfv12}.
Following \cite{7}, the parasitic growth parameter can be calculated by the formula,
\begin{equation}\label{NA9}
    \mu_i =  (\xi_i)^{-1} w_i^*BUw_i,   
\end{equation}
\noindent where $\xi_i$ is the $i$-th eigenvalue of V and $ w_i $ is the corresponding left eigen vector with $\xi_1=1$ and $\xi_i\not= 1$  for $2\leq i\leq r$.

\section{Projection technique for general linear methods}
\noindent Let $y' = f(y(x))$ denote a differential equation on a manifold $\textbf{M}$, with $g(y)$ as an invariant, such that,
\begin{align*}
    \textbf{M}= \{y ; g(y) = 0\}.
\end{align*}
The exact solution stays on the manifold $M$,
\begin{align} \label{ABC22}
 y_{0} \in \textbf{M}    \implies y(x)\in \textbf{M} \quad \forall x.  
\end{align}
We want our numerical solution by the general linear method to stay on the manifold $M$. For this purpose we employ standard projection technique for general linear methods \cite{7}. 

\noindent Let $y^{[m]} \in \textbf{M}$. An application of one step of the GLM yields $ \tilde {y}^{[m+1]}\notin \textbf{M}$. Project the value $ \tilde {y}^{[m+1]}_{1}$ onto the manifold $\textbf{M}$ to obtain $y^{[m+1]}_{1} \in \textbf{M}$ such that,
    \begin{equation*} \label{ABCD1}
     y^{[m+1]}_{1}= \tilde {y}^{[m+1]}_{1} +  \frac{H(y_{0})-H( \tilde {y}^{[m+1]}_{1})}{<\nabla H(\tilde{y}^{[m+1]}_{1}),\nabla H(\tilde{y}^{[m+1]}_{1}) >} \nabla H(\tilde{y}^{[m+1]}_{1})   
     \end{equation*}
\noindent where $\nabla H(y)$ is the gradient of $H(y)$. The important observation is that the projection method is applied on the first output value  $y^{[m+1]}_{1}$ only.
\section{ Variational integrators as general linear methods}

\noindent In order to shed light on variational integrator for degenerate Lagrangian system (\ref{ABC2}) expressed as a general linear method we consider non-linear pendulum whose Lagrangian is degenerate, 
\begin{gather}\label{dlp}
  L(q, \dot{q}) = \begin{bmatrix} q^{[2]} & 0  \end{bmatrix} \begin{bmatrix} \dot{q}^{[1]}\\ \dot{q}^{[2]}    \end{bmatrix} +cos(q^{[1]}) - \frac{(q^{[2]})^2}{2}.
 \end{gather}
Comparing \eqref{dlp} with  \eqref{dl} we get,
\begin{gather*}
     \alpha (q)= \begin{bmatrix}
q^{[2]} & 0 
\end{bmatrix},\hspace{0.5in} H(q)= \frac{(q^{[2]})^2}{2} - cos(q^{[1]}).
\end{gather*}
Consequently,

\begin{gather*}
 \nabla \alpha (q) = \begin{bmatrix} 0 & 0 \\ 1 & 0
 \end{bmatrix}, \hspace{0.6in}  \nabla H(q)=  \begin{bmatrix}
sin(q^{[1]}) \\ q^{[2]}
\end{bmatrix}.
 \end{gather*}
By inserting these values in (\ref{ABC2}), we get,
\begin{gather*}
      \begin{bmatrix}
q_{m+1} ^{[1]} \\ q_{m+1} ^{[2]}
\end{bmatrix} =  \begin{bmatrix}
q_{m-1}^{[1]} \\ q_{m-1}^{[2]} 
\end{bmatrix} + 2h \begin{bmatrix} q_{m}^{[2]} \\ -sin(q_{m}^{[1]})
\end{bmatrix},
\end{gather*}

\begin{align}\label{vip}
 \implies y_{m+1}=y_{m-1}+2hf.  
\end{align}
The equation \eqref{vip} is the required variational integrator for degenerate Lagrangian system representing non-linear pendulum. Evidently, equation \eqref{vip} is a multistep method which can be written as general linear method \eqref{Adams_Moulton} as,
 \begin{align}\label{glmp}
 \begin{bmatrix}
    \begin{array}{c|c}
    \boldsymbol{\textbf{A}} & \boldsymbol{\textbf{U}} \\ \hline \ \boldsymbol{\textbf{B}} & \boldsymbol{\textbf{V}}
\end{array}
\end{bmatrix} =
      \begin{bmatrix}
    \begin{array}{c|cccc}
       0 & 0 & 1 & 2 & 0 \\
       \hline
     0 & 0 & 1 & 2 & 0 \\
      0 & 1 & 0 & 0 & 0 \\
       1 & 0 & 0 & 0 & 0 \\
       1 & 0 & 0 & 1 & 0 \\
    \end{array}
    \end{bmatrix}.
 \end{align}
 The parasitic growth parameters of \eqref{glmp} by using \eqref{NA9} is computed as $\mu = -1.667$.
 \subsection{Starting Algorithm}   
\noindent To find the value of $q_{-1}$, we use the position momentum form \cite{NA19,11},
\begin{align}
    p_m&=-D_1 L_d(q_m,q_{m+1}),
\\  \label{A}  p_{m+1}&= D_2 L_d(q_m,q_{m+1}).
\end{align}
To obtain a relation between $q_{-1}$, $q_{0}$ and $p_{0}$, use  the equation (\ref{A}) as,
\begin{equation*}
    p_{0}=D_{2}L_{d}(q_{-1},q_{0}),
\end{equation*}
 but,
 
 \begin{equation*}
 p_{0}=\alpha(q_{0}),     
 \end{equation*}

\begin{equation*}
\implies    \alpha(q_{0})= D_{2}L_{d}(q_{-1},q_{0}).
\end{equation*}
     
\section{Numerical experiment}
\begin{figure} [H] 
\centering
 \includegraphics[width=110mm]{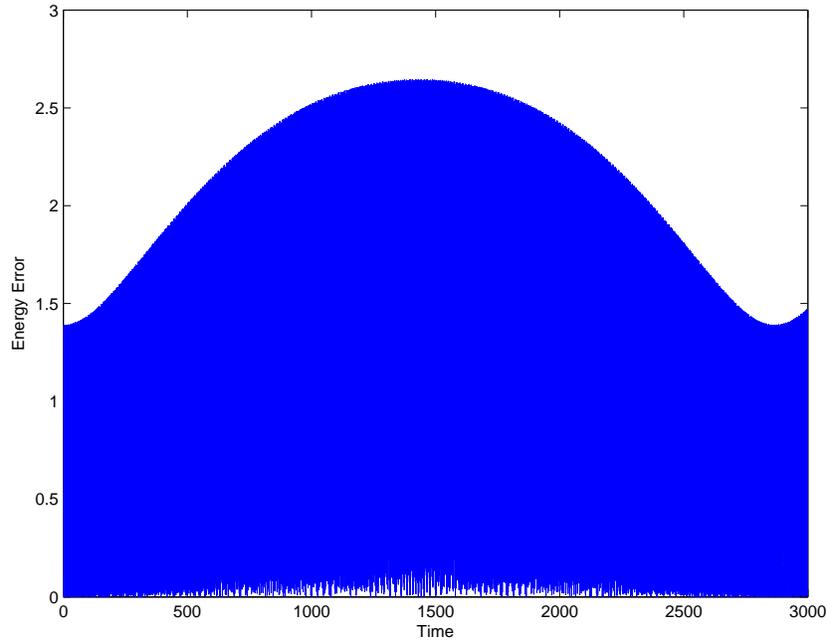}
    \caption{Energy error in non-linear pendulum by using variational integrator without projection.}
    \label{fig:6}
\end{figure}

\begin{figure} [H]
    \centering
    \includegraphics[width=110mm]{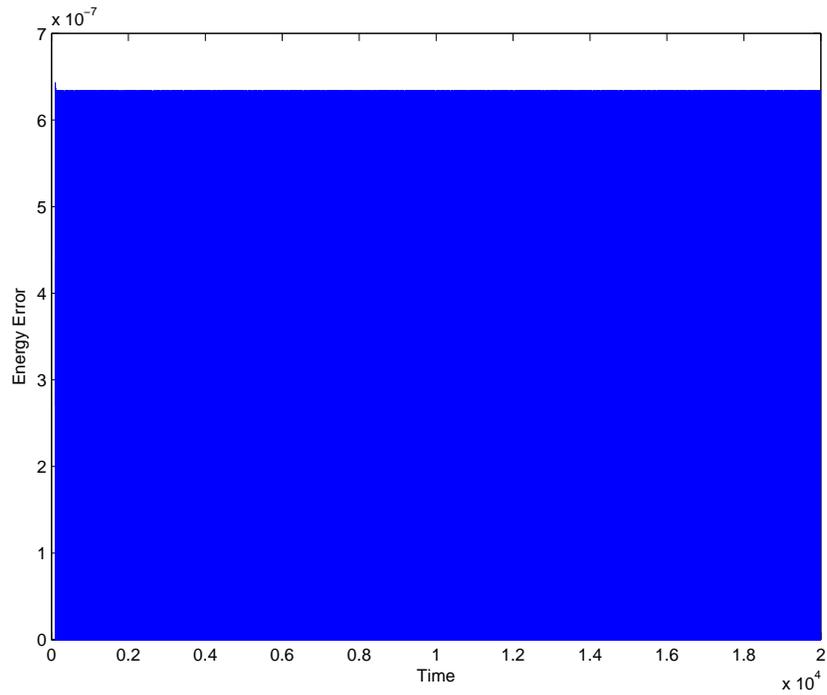}
    \caption{Energy error in non-linear pendulum by using variational integrator with projection.}
    \label{fig:7}
    \end{figure}

 \noindent An application of \eqref{glmp} with initial condition $q_{0}=(2.3,0)$ and step-size $h=0.1$ yields energy error in the pendulum and is shown in Figure \ref{fig:6}.
\noindent Figure \ref{fig:6} shows that the variational integrator \eqref{glmp} does not conserve the energy. We have calculated the absolute error as follows,
\[ Error=abs(H_e-H_n),\]
\noindent where ``$H_e$" is the exact energy at initial point and ``$H_n$ " is the approximate energy calculated at all numerical values. We then apply the projection technique on GLM \eqref{glmp} and calculate the energy error again as shown in Figure \ref{fig:7}.
\noindent Figure \ref{fig:7}  shows that the variational integrator for degenerate Lagrangin of non-linear pendulum as GLM \eqref{glmp} with projection technique is preserving the energy very well.

\end{document}